\documentclass[a4paper, oneside, 11pt]{article}
\usepackage[english]{babel}
\usepackage{amsfonts}
\usepackage{amssymb}
\usepackage{graphics}
\usepackage{amsrefs}
\usepackage{latexsym}
\begin{document}

\title{{\bf{\Large{Bayesian nonparametric estimation of Simpson's evenness index under $\alpha-$Gibbs priors.}}}\footnote{{\it AMS (2000) subject classification}. Primary: 60G58. Secondary: 60G09.}}
\author{\textsc {Annalisa Cerquetti}\footnote{Corresponding author, SAPIENZA University of Rome, Via del Castro Laurenziano, 9, 00161 Rome, Italy. E-mail: {\tt annalisa.cerquetti@gmail.com}}\\
\it{\small Department of Methods and Models for Economics, Territory and Finance}\\
  \it{\small Sapienza University of Rome, Italy }}
\newtheorem{teo}{Theorem}
\date{\today}
\maketitle{}

\begin{abstract}
A Bayesian nonparametric approach to the study of species diversity  based on choosing a random discrete distribution as a prior model for the unknown relative abundances of species has been recently introduced in Lijoi et al. (2007, 2008). Explicit posterior predictive estimation of {\it species richness} has been obtained under priors belonging to the $\alpha$-Gibbs class (Gnedin \& Pitman, 2006). Here we focus on posterior estimation of {\it species evenness} which accounts for diversity in terms of the proximity to the situation of uniform distribution of the population into different species. We focus on Simpson's index and  provide a Bayesian estimator under quadratic loss function, with its variance, under some specific $\alpha-$Gibbs priors.
\end{abstract}

\section*{Introduction}

Species diversity is made of two components: species {\it richness}, the number of different species belonging to a population, and species {\it evenness}, the distance of the actual relative abundances from a situation of uniform distribution of the population into different species. In other words a population of species is the more {\it diverse} the more is {\it rich} and {\it uniform}. The measurement of diversity of populations when individuals are classified into groups has a long history, dating back to Simpson's (1949) and Fisher's (1943) seminal papers. Since then ecological literature has produced a huge amount of results and a variety of indexes to measure both aspects of diversity. With respect to {\it richness} nonparametric estimation relevant recent approaches includes Chao \& Lee (1992), Shen et al. (2003) and Chao \& Bunge (2002), while a Bayesian approach is in Hill (1979) and Boender \& Rinnooy Kan (1987). With respect to {\it evenness} nonparametric estimation, results mostly concentrate on Simpson's (1949) index and Shannon-Wiener's index (Pielou, 1975). (See e.g. Lloyd \& Ghelardi, 1964, and Chao \& Shen, 2003). A Bayesian nonparametric estimation of Shannon's index under Fisher's (1943) model is in  Gill \& Joanes (1979). \\

Recently a Bayesian nonparametric approach to species diversity has been introduced in Lijoi et al. (2007, 2008) by specifying a {\it prior} on the space of unknown relative abundances and deriving posterior results given the information contained in the multiplicities of the first $k$ species observed in a basic $n$-sample. Posterior predictive estimation of species {\it richness}  have been obtained under {\it priors} belonging to random discrete distributions inducing infinite exchangeable partitions in Gibbs form of type $\alpha$, as devised in Gnedin \& Pitman (2006).  Asymptotic results for the behaviour of a proper normalization of the number of new species when the size of the additional sample tends to infinity have been provided under the same approach (Favaro et al. 2009, 2011; Cerquetti, 2011).

Here, under the same setting, we propose a Bayesian nonparametric posterior estimation of the {\it evenness} of a population of species by focusing on {\it Simpson's index} of species diversity (Simpson, 1949) which actually measures evenness by the chance to observe two elements picked at random belonging to different species. We derive an explicit posterior estimator under quadratic loss function with its variance under {\it two parameter $(\alpha, \theta)$ Poisson-Dirichlet model} (Pitman, 1995, Pitman and Yor, 1997) and consequently derive explicit formulas for three particular cases this model contains: Dirichlet-Ewens (1972) sampling model ($\alpha=0$), Stable model ($\theta =0$) and Fisher's model (1943), ($\theta=\xi|\alpha|$ for $\alpha <0$ and $\xi \in \mathbb{N}^+$).  The organization of the paper is as follows. In Section 1 we recall some preliminaries about $\alpha-$Gibbs exchangeable partitions and Simpson's index. In Section 2 we locate the {\it evenness} parameter by providing prior estimator of Simpson's index under  some relevant $\alpha-$Gibbs partition  models. In Section 3 we obtain the general Bayesian posterior nonparametric estimator under the two-parameter $(\alpha, \theta)$ model and its posterior variance, and derive the explicit posterior results for the particular cases treated in Section 2. Finally in Section 4  we suggest a feasible generalization of the results contained in the present paper.

\section{Preliminaries}

Recall that the exchangeable partition probability function describing the random behaviour of a partition of the natural integers obtained by sampling from a random discrete distribution $(P_i)_{i \geq 1}$ belonging to the $\alpha$-Gibbs class for $\alpha \in (-\infty, 1)$ (Gnedin and Pitman, 2006) is characterized by the product form
\begin{equation}
\label{gibbseppf}
p(n_1, \dots, n_k)=V_{n,k} \prod_{j=1}^k (1 -\alpha)_{n_j -1}
\end{equation}
where  $V_{n,k}$ are coefficients satisfying the recursive relation $V_{n,k}= (n -k\alpha) V_{n+1, k}+ V_{n+1, k+1}$, and $(x)_y$ is the usual notation for {\it rising} factorials $(x)_y=x(x+1)\cdots(x+y-1)$. Three convex subclasses of $\alpha$-Gibbs partitions have been identified by Gnedin \& Pitman by deriving extreme partitions corresponding to different values of $\alpha$, namely: mixtures of Fisher's models for $\alpha <0$, mixtures of Ewens-Dirichlet ($\theta$) models for $\alpha=0$, and mixtures of  conditional Poisson-Kingman models driven by the $\alpha$ stable subordinator (Pitman, 2003) for $\alpha \in (0,1)$. 

The Bayesian nonparametric approach to species sampling problems under $\alpha$-Gibbs priors has been introduced in Lijoi et al. (2007, 2008) under the following observational set up. Given $X_1, \dots, X_n$ a sample of observations (labels) from an infinite population of species with unknown relative abundances $(P_i)_{i \geq 1}$ which has provided  $k$ different species with multiplicities in order of appearance $(n_1, \dots, n_k)$, let $X_{n+1}, \dots, X_{n+m}$ be a future sample of $m$ additional observations. Then interest lies in posterior (conditional) inference on some characteristic of the unknown population of abundances or in predictive inference on some characteristic of a future sample. Lijoi et al. (2007, 2008) provide the predictive distribution  under a general $\alpha-$Gibbs prior, for various quantities of interest related to the additional sample. In particular for $K_m$ the number of {\it new species} they obtain the following result
\begin{equation}
\label{richness}
\mathbb{P}_{\alpha, V_{n,k}}(K_m=k^*|n_1, \dots, n_k)=\frac{V_{n+m, k+k^*}}{V_{n,k}} S_{m, k^*}^{-1, -\alpha, -(n-k\alpha)}
\end{equation}
where $k^*$ takes values in $0, \dots, m$ and $S_{m, k^*}^{-1, -\alpha, -(n -k\alpha)}$ are non central generalized Stirling numbers. 
Now the most studied and tractable $\alpha-$Gibbs model is the {\it two-parameter Poisson-Dirichlet} $(\alpha, \theta)$ model, for $\alpha \in (0,1)$ and $\theta>-\alpha$ or $\alpha <0$ and $\theta=|\alpha|\xi$ for $\xi=1, 2,\dots $ introduced in Pitman (1995) and largely studied in Pitman and Yor (1997). It arises by mixing Poisson-Kingman models driven by the Stable subordinator (Pitman, 2003) by a polynomial tilting of the stable density and its EPPF is given by
\begin{equation}
\label{twopar}
p_{\alpha, \theta} (n_1, \dots, n_k)= \frac{(\theta + \alpha)_{k -1 \uparrow \alpha}}{(\theta +1)_{n-1}} \prod_{j=1}^k (1- \alpha)_{n_j-1},
\end{equation}
where $(x)_{y \uparrow \alpha}$ stands for generalized rising factorials $(x)_{y \uparrow \alpha}=(x)(x+\alpha) \cdots(x +(y-1)\alpha )$. A huge amount of results is available for this model, mostly due to J. Pitman, and we will heavily rely on them in what follows. The posterior predictive distribution for the number of new species in an additional $m$-sample easily follows from (\ref{richness}) by some elementary combinatorial calculus and corresponds to
$$
\mathbb{P}_{\alpha, \theta}(K_m=k^*|n_1, \dots, n_k)=\frac{(\theta +k\alpha)_{k^* \uparrow \alpha}}{(\theta +n)_m} S_{m,k}^{-1, -\alpha, -(n -k\alpha)}.
$$
Additional results for an explicit Bayesian estimator of predictive species richness and on the asymptotic behaviour of a proper normalization of $K_m$  under two-parameter Poisson-Dirichlet priors are in Favaro et al. (2009), (see Cerquetti, 2011, for a generalization of the asymptotic result to general $\alpha-$Gibbs partitions for  $\alpha \in (0, 1)$).\\\\
In his seminal paper on {\it Nature} in 1949 Simpson introduces an index of diversity (''a measure of the concentration of the classification'') for an infinite population such that each individual belongs to one of $k$ different groups, for $P_1, \dots, P_k$ the proportions of individuals in variuos groups. It is defined as  
$$
H_S= 1 - \sum_{j=1}^k P_j^2
$$
namely the {\it chance that two observations picked at random from the population belong to different species}. The index can take any values between $0$ and $1-1/k$, the former representing the smallest evenness ({\it maximum concentration}) and the latter the situation of maximum diversity ({\it minimum concentration}), hence it depends both on the richness of the population like on the evenness. Here dealing with priors on discrete distributions on an infinite number of classes we work with the infinite species version of Simpson's index which vares in $(0, 1)$ hence becomes just an index of {\it evenness}. 

\section {Prior evenness in $\alpha-$Gibbs partition models}

Apart from Fisher's model and mixtures of Fisher's models (Gibbs class for $\alpha <0$), choosing a nonparametric prior in the $\alpha$-Gibbs class, actually implies assuming the number of different species in the population under study to be unknown and theoretically {\it infinite}.  It follows that if we want to introduce a prior knowledge on the diversity of the population we cannot act on any kind of {\it richness} parameter. Nevertheless for each $\alpha$-Gibbs prior it is always possible to identify an {\it evenness} parameter.  Here we show how by studying  the prior mean of Simpson's index. \\\\
{\bf Proposition 1.} [Prior mean and variance of Simpson index under $\alpha-$Gibbs priors] {\it Let $(X_1, \dots, X_n)$ be a sample from a population of species whose relative abundances are distributed according to a general $\alpha-$Gibbs partition model,  then for $S_2=\sum_{j=1}^\infty P_j^2$, the sequence of prior moments of $1-H_S$ is given by
\begin{equation}
\label{momsim}
\mathbb{E}_{\alpha, V_{n,k}}[S_2^\xi]=\sum_{j=1}^\xi \frac{1}{j!} V_{2\xi, j} \sum_{\xi_1, \dots, \xi_j} \frac{\xi!}{\xi_1! \cdots, \xi_j!} \prod_{i=1}^j (1 -\alpha)_{2\xi_i-1},
\end{equation}
hence a prior estimate of Simpson's index corresponds to
\begin{equation}
\label{pr_mean}
\mathbb{E}_{\alpha, V_{n,k}}(H_S)=1-\mathbb{E}(\sum_{j=1}^\infty \tilde P_j^2)=1-\mathbb{E}(S_2)=1-V_{2,1}(1-\alpha)
\end{equation}
and its prior variance is given by}
\begin{equation}
\label{pr_var}
Var(H_S)=Var(S_2)=\mathbb{E}(S_2^2)-[\mathbb{E}(S_2)]^2=V_{4,1}(1-\alpha)_3 + V_{4,2}(1-\alpha)^2- [V_{2,1}(1-\alpha)]^2.
\end{equation}
{\it Proof}. Let $(P_i)_{i \geq 1}$ be the sequence of ranked atoms of a random discrete distribution, and $\tilde{P}_j$ the random size of the $j$th atom discovered in the process of random sampling, or equivalently the asymptotic frequency of the $j$th class when the blocks of the partition generated are put in order of their least element.  Now for the random variable 
$$
S_m:=\sum_{i=1}^\infty P_i^m= \sum_{j=1}^\infty \tilde{P}_j^m,
$$
where it is still assumed that $S_1=1$ almost surely, Pitman (2003) provides the following general expression for the $\xi$th  moment 
\begin{equation}
\label{moments}
\mathbb{E}[S_m^\xi]= \sum_{j=1}^\xi \frac{1}{j!} \sum_{\xi_1, \dots, \xi_j} \frac{\xi!}{\xi_1! \cdots \xi_j!} p(m\xi_1, \dots, m\xi_j), 
\end{equation}
where the second sum is over all sequences of $j$ positive integers $(\xi_1, \dots,\xi_j)$ with $\xi_1+\dots+\xi_j=\xi$. This implies the EPPF induced by sampling from a random discrete distribution directly determines the positive integers moments of the power sums $S_m$, hence the distribution of $S_m$ for each $m$. For $m=2$, by substituting the Gibbs form (\ref{gibbseppf}) in (\ref{moments}), (\ref{momsim}) follows. For $\xi=1$ and $\xi=2$  (\ref{pr_mean}) and (\ref{pr_var}) can be easily obtained. \hspace{11.8cm}$\square$ \\\\
{\bf Example 2.} [{\it Two parameter Poisson-Dirichlet model}]. For $(P_j)_{j \geq 1}$ having two parameter Poisson-Dirichlet distribution $(\alpha, \theta)$, for $\alpha \in (0,1)$ and $\theta > -\alpha$ then its size-biased permutation is characterized by (Pitman, 1996)
$$
\tilde{P}_j=W_i \prod_{i=1}^{j-1} (1 -W_i)
$$ 
for $W_i$ independent $Beta(1- \alpha, \theta +j\alpha)$. The prior mean of Simpson's index hence corresponds to 
$$
\mathbb{E}_{\alpha, \theta}(H_s)=\frac{\theta+\alpha}{1 +\theta},
$$
which is the probability to observe a {\it  new} species in the classical sequential Chinese restaurant construction of the partition given the first observation.  It follows both large values of $\theta$ and $\alpha$ may increase the prior belief on the evenness of the population. As for the prior variance
$$
Var_{(\alpha, \theta)}(H_s)= \frac{(1 -\alpha)_3+ (\theta +\alpha)(1-\alpha)^2}{(\theta +1)_3} - \frac{(1-\alpha)^2}{(\theta +1)^2}.
$$
{\bf Example 3.} [{\it Fisher's model}]  The simpler and older species sampling model was introduced in Fisher  et al. (1943), and it is characterized by a finite number $\xi$ of species whose unknown proportions  $(P_1, \dots, P_\xi)$ have a symmetric Dirichlet distribution with parameter $(|\alpha|, \dots, |\alpha|)$, for $\alpha \in (-\infty, 0)$, arising by $(G_1/G, \dots, G_\xi/G)$ for $G_i$ independent $Ga(|\alpha|, 1)$ r.v.s with $G=\sum_i G_i$. To identify the {\it evenness} parameter is enough to notice that the marginal of $P_j$ are $Beta(|\alpha|, \xi|\alpha|)$ distributed, so that 
$$
\mathbb{E}_{\alpha, \xi}(H_s)=1-\frac{1+\alpha}{1+\xi|\alpha|}=\frac{\xi|\alpha|-\alpha}{1 +\xi|\alpha|}.
$$
Hence large values of $|\alpha|$ means $(P_i)$ fairly uniform are expected, while small values of $|\alpha|$ correspond to widely different $P_i$ expected. This shows $|\alpha|$ is indeed an evenness parameter. Its explicit prior variance is as follows
$$
Var_{\alpha, \xi} (H_s)= \frac{(1 + \alpha)_3+ (|\alpha|\xi - |\alpha|)(1 +\alpha)^2}{(|\alpha|\xi+1)_3} - \frac{(1+\alpha)^2}{(|\alpha|\xi +1)^2}.
$$ 
{\bf Example 4.} [{\it Dirichlet-Ewens model}]. For $\alpha \rightarrow 0$, $\xi \rightarrow  \infty$ and $|\alpha|\xi\rightarrow\theta$ then Fisher's model converges to the random atoms of the  Dirichlet probability measure  (Ferguson, 1973), whose ranked atoms are known to have  $(P_j)\sim PD(\theta)$ Poisson-Dirichlet distribution. Their size biased permutation has GEM distribution
$$
\tilde{P_j}=W_j \prod_{i=1}^{j-1} (1 -W_i)
$$
for $W_i$ iid $\sim Beta(1, \theta)$, hence 
$$
\mathbb{E}_{\theta}(H_s)= \frac{\theta}{1+\theta}.
$$
This shows that $\theta$ is exactly a parameter of evenness and that choosing a Dirichlet prior over the unknown relative abundances allows to explicitly introduce a prior hypothesis on the evenness of the population. The prior variance of $H_S$ corresponds to 
$$
Var_{\theta}(H_s)= \frac{6+ \theta}{(\theta+1)_3} - \frac{1}{(\theta +1)^2}.
$$
{\bf Remark 5.} Notice that mixing over $\theta$, as from Gnedin \& Pitman (2006), with a general mixing distribution over $(0, \infty)$ produces the class of exchangeable Gibbs partitions of type $\alpha=0$. In terms of random probability measure this corresponds to Antoniak's  (1974) mixtures of Dirichlet processes. Despite it has received less attention in the Bayesian treatment of species sampling problems, we notice this class actually entails the possibility to explicitly introduce a prior on the evenness parameter of the Dirichlet partition model. \\\\
{\bf Example 6.} [{\it $\alpha$-Stable model}] From the results for the Poisson-Dirichlet $(\alpha, \theta)$ model,  prior mean and variance of the Simpson's index for random discrete distribution obtain by normalization of the ranked lengths of the ranked points of a Poisson process with mean intensity the L\'evy density of the $\alpha-$Stable density for $\alpha \in (0,1)$, easily follow by letting $\theta=0$, since $PK(\rho_\alpha)=PD(\alpha, 0)$, hence 
$$
\mathbb{E}_{\alpha}(H_S)={\alpha}
$$
and
$$
Var_{\alpha}(H_s)= \frac{(1 -\alpha)_3 + (\alpha)(1 -\alpha)^2}{6} - (1-\alpha)^2= \frac{2\alpha(1-\alpha)}{6}.
$$
It follows even in this case the value of $\alpha$ is a direct measure of the evenness of the population and may be chosen accordingly to prior information in a Bayesian nonparametric perspective.

\section{Posterior estimation of species evenness under PD$(\alpha, \theta)$ priors}

To obtain a Bayesian nonparametric estimator of Simpson's index of evenness we just need the posterior distribution of the size-biased ordered atoms of the random discrete distribution chosen as the prior model. An explicit result for the two parameter $(\alpha, \theta)$ model is in Pitman (1996). Notice that all the Examples treated in the previous section arise as particular cases of this model hence specific posterior estimators and their variance may be obtained by the general result for this class.  \\\\
{\bf Proposition 7.} {\it Let $(X_i)_{i \geq 1}$ be a population of exchangeable species labels driven by a two-parameter Poisson-Dirichlet prior. Let $(n_1, \dots, n_k)$ be the multiplicities of the first $k$ species observed in a basic $n$-sample. Then the Bayesian nonparametric estimate under quadratic loss function of Simpson's index of evenness $H_S$  is given by}
\begin{equation}
\label{post_mean}
\mathbb{E}_{\alpha, \theta} (H_s|(n_1, \dots, n_k))= 1 -  \frac{\sum_{j=1}^k(n_j -\alpha)_{2}}{(\theta +n)_2}+\frac{(\theta +k \alpha)(1-\alpha)}{(\theta +n)_2}. 
\end{equation}
{\it Proof.} Recall from Pitman (1996) that if $(X_n)$ is a sample from a random discrete distribution $P$ with  ranked atoms having $PD(\alpha, \theta)$ distribution, then, conditionally given $(n_1, \dots, n_k)$, the random partition of $n$ induced by the first $k$ different values $\tilde{X}_1, \dots, \tilde{X}_k$ of $(X_n)$, the posterior random discrete distribution is given by
\begin{equation}
\label{postmod}
P_n(\cdot)=\sum_{j=1}^k \tilde{P}_{j, {\bf n}} \delta_{\tilde{X}_j}(\cdot) + \tilde{R}_k P_k(\cdot),
\end{equation}
where $(\tilde{P}_{1, {\bf n}}, \dots, \tilde{P}_{k,{\bf n}} \tilde{R}_k)$ has Dirichlet $(n_1 -\alpha, \dots, n_k -\alpha, \theta +k\alpha)$ distribution, {\it independently} of the random discrete distribution $P_k(\cdot)$ which has ranked atoms  $(Q_i)_{i \geq 1}$ having $PD(\alpha, \theta +k\alpha)$ distribution. Now, by the independence
$$
\mathbb{E}\left(\sum_{j=1}^\infty P_{j,{\bf n}}^{*2}\right)= \mathbb{E} (\sum_{j=1}^k \tilde {P}_{j, {\bf n}}^2) + \mathbb{E}(\tilde{R}_k^2) \mathbb{E} \left(\sum_{i=1}^{\infty} {Q}_{i, \theta +k \alpha}^2 \right)
$$
and recalling that if $X \sim Beta(a, b)$ then $E(X^k)=\frac{(a)_k}{(a+b)_k}$,  since $\tilde{P}_{j,{\bf n}} \sim Beta(n_j -\alpha, \theta +n - n_j+\alpha)$, and $\tilde{R}_k \sim Beta(\theta +k\alpha, n -k\alpha)$ then
$$
\mathbb{E} (\sum_{j=1}^k \tilde {P}_{j, {\bf n}}^2) =\sum_{j=1}^k \mathbb{E}\tilde{P}_{j,{\bf n}}^2=\frac{\sum_{j=1}^k  (n_j -\alpha)_{2}}{(\theta +n)_2}.
$$
Additionally, by the structural distribution of $PD(\alpha, \theta +k\alpha)$, $\tilde{Q}_1 \sim Be(1-\alpha, \theta +k\alpha +\alpha)$ hence
$$
\mathbb{E}(\tilde{R}_k^2)\mathbb{E} \left(\sum_{i=1}^{\infty} {Q}_{i, \theta +k \alpha}^2 \right)= \frac{(\theta +k\alpha)_2}{(\theta +n)_2} \frac{1-\alpha}{\theta +k\alpha +1}
$$
and the result is proved. $\hspace {9.2cm}$ $\square$\\\\
It is even possible to obtain the posterior variance of the Bayesian estimator of the Simpson's index under two parameter Poisson-Dirichlet model. \\\\
{\bf Proposition 8.} {\it Let $(X_i)_{i \geq 1}$ be a population  of exchangeable species labels driven by a two-parameter Poisson-Dirichlet prior. Let $(n_1, \dots, n_k)$ be the multiplicities of the first $k$ species observed in a basic $n$-sample, then the Bayesian nonparametric estimator of Proposition 7. is characterized by the following posterior variance} 
\begin{equation}
\label{post_var}
Var_{\alpha, \theta}(H_s|(n_1, \dots, n_k))=\frac{\sum_{j=1}^k (n_j -\alpha)_4}{(\theta +n)_4} + \frac{2 \sum_{i \neq j}(n_j -\alpha)_2(n_i -\alpha)_2}{(\theta +n)_4}+
$$
$$
+\frac{(\theta +k\alpha)}{(\theta +n)_4}\frac{[(1-\alpha)_3 +(\theta +k\alpha +\alpha) (1- \alpha)^2]}{}+$$
$$+ \frac{2\sum_{j=1}^k (n_j -\alpha)_2(\theta + k\alpha)}{(\theta +n)_4} \frac{(1 -\alpha)}{} - \left[ \frac{\sum_{j=1}^k (n_j -\alpha)_2 + (\theta +k\alpha)(1-\alpha)} {(\theta +n)_2}\right]^2.
\end{equation}
{\it Proof:} It is enough to obtain the posterior second moment of $1-H_s$. We can always write
$$
\mathbb{E}_{\alpha, \theta} \left(\sum_{j=1}^{\infty} \tilde{P}^{*2}_{j,{\bf n}} \right)^2= \mathbb{E}_{\alpha, \theta} \left(\sum_{j=1}^k \tilde{P}^2_{j,{\bf n}} + R_k^2 \sum_{i=1}^{\infty} \tilde{Q}_i^2 \right)^2
$$
which reduces to 
$$
=\sum_{j=1}^k \mathbb{E}(\tilde{P}^4_{j, {\bf n}})+ 2 \sum_{i \neq j}\mathbb{E}(\tilde{P}^2_{i, {\bf n}}\tilde{P}^2_{j, {\bf n}})+\mathbb{E}(R_k)^4\mathbb{E}(\sum_{i=1}^\infty \tilde{Q}^2_{i})^2 
+ 2 \sum_{j=1}^k \mathbb{E}(\tilde{P}^2_{j, {\bf n}} R_k^2) \mathbb{E}(\sum_{i=1}^{\infty} \tilde{Q}_i^2).
$$
By known expressions for mixed moments of Dirichlet vectors, and since  from Pitman's formula recalled in Section 1, 
$$
\mathbb{E}(\sum_{i=1}^\infty \tilde{Q}_i^2)^2=\mathbb{E}(S_{2, {\theta +k\alpha}})^2= V_{4,1}(1-\alpha)_3 + V_{4,2}(1-\alpha)^2= \frac{(1-\alpha)_3 +(\theta +k\alpha +\alpha) (1- \alpha)^2}{(\theta + k\alpha +1)_{3}},
$$
it follows
$$
\mathbb{E}_{\alpha, \theta} \left(\sum_{j=1}^{\infty} \tilde{P}^{*2}_{j,{\bf n}} \right)^2=\frac{\sum_{j=1}^k (n_j -\alpha)_4}{(\theta +n)_4} + \frac{2 \sum_{i \neq j}(n_j -\alpha)_2(n_i -\alpha)_2}{(\theta +n)_4}+
$$
$$
+\frac{(\theta +k\alpha)}{(\theta +n)_4}\frac{[(1-\alpha)_3 +(\theta +k\alpha +\alpha) (1- \alpha)^2]}{}+
\frac{2\sum_{j=1}^k (n_j -\alpha)_2(\theta + k\alpha)}{(\theta +n)_4} \frac{(1 -\alpha)}{}.
$$
\hspace{13.7cm}$\square$\\\
{\bf Example 9.} [{\it Fisher's model (continued)}] The Bayesian nonparametric estimator for Simpons's index of evenness under $(\alpha, \xi|\alpha|)$ Fisher's model is given by
$$
\mathbb{E}_{\xi, \alpha}(H_S|(n_1, \dots, n_k))= 1 - \frac{\sum_{j=1}^k (n_j -\alpha)_2 + (|\alpha|\xi -k\alpha)(1 +\alpha)}{(|\alpha|\xi +n)_2}
$$
and its posterior variance corresponds to 
$$
Var(H_S|(n_1, \dots, n_k))= \frac{\sum_{j=1}^k (n_j +\alpha)_4}{(|\alpha|\xi +n)_4} + \frac{2 \sum_{i\neq j} (n_j +\alpha)_2(n_i +\alpha)_2}{(|\alpha|\xi +n)_4}+
$$
$$
+ \frac{(|\alpha|\xi -k\alpha)[(1 +\alpha)_3 +(|\alpha|\xi -k \alpha -\alpha)(1+\alpha)^2]}{(|\alpha|\xi +n)_4} + \frac{2 \sum_{j=1}^k(n_j +\alpha)_2 (|\alpha|\xi - k\alpha)(1 +\alpha)}{(|\alpha|\xi +n)_4} - $$
$$-\left[ \frac{\sum_{j=1}^k (n_j +\alpha)_2 + (|\alpha|\xi -k\alpha)(1+\alpha)}{(|\alpha|\xi +n)_2}\right]^2
$$
{\bf Example 10} [{\it Dirichlet-Ewens model (continued)}] The Bayesian nonparametric estimator for Simpson's index arises specializing (\ref{post_mean}) for $\alpha=0$ hence 
$$
\mathbb{E}_{\theta}(H_s|(n_1, \dots, n_k))= 1 - \frac{\sum_{j=1}^k (n_j)_2 + \theta }{(\theta +n)_2}.
$$
Similarly for the posterior variance of $H_S$ we obtain
$$
Var_{\theta}(H_s|(n_1, \dots, n_k))= \frac{\sum_{j=1}^k (n_j)_4 + 2 \sum_{i \neq j} (n_j)_2(n_i)_2 +\theta[(1)_3 + \theta ] + 2 \theta \sum_{j=1}^k (n_j)_2}{(\theta+n)_4} - 
$$
$$-\left(\frac{(\sum_{j=1}^k (n_j)_2 + \theta }{(\theta +n)_2}\right)^2.
$$\\\\
{\bf Example 11.} [{\it $\alpha-$Stable model (continued)}] Posterior mean and variance of Simpson's index under $\alpha-$Stable model for $\alpha \in (0,1)$ arise by the general formulas for $\theta=0$ hence
$$
\mathbb{E}_{\alpha}(H_s| n_1, \dots, n_k)= 1- \frac{\sum_{j=1}^k (n_j -\alpha)_2}{(n)_2}+\frac{k\alpha(1 -\alpha)}{(n)_2}
$$
and
$$
Var(H_S|{n_1, \dots, n_k})= \frac{\sum_{j=1}^k (n_j -\alpha)_4}{(n)_4} +
$$
$$+\frac{2 \sum_{i \neq j} (n_j -\alpha)_2 (n_i -\alpha)_2}{(n)_4} + \frac{k\alpha[(1 -\alpha)_3+(k\alpha +\alpha)(1- \alpha)^2]}{(n)_4}+ $$
$$+\frac{2 \sum_{j=1}^k (n_j -\alpha)_2 k\alpha(1- \alpha)}{(n)_4} - \left[\frac{\sum_{j=1}^k (n_j -\alpha)_2 +k\alpha(1- \alpha)}{(n)_2}\right]^2.
$$

\section {Future directions}

A generalization of  Pitman's posterior result for the two parameter model (\ref{postmod}) for random probability measure obtained by normalization of measures with independent increments is in James et al. (2009).  Additionally an explicit stick breaking construction for the size biased atoms of the normalized Inverse Gaussian partition model, which belongs to the $\alpha$ Gibbs class for $\alpha=1/2$ and mixing distribution the exponentially tilted version of the $1/2-$Stable distribution (cfr. Pitman, 2003), has been recently obtained in Favaro et al. (2012).  Those results suggest it may be worth investigating the possibility to obtained an explicit Bayesian nonparametric estimator for Simpson's index of evenness under this prior model for relative abundances. This will be the subject of future investigations.

\section*{References}
\newcommand{\bibu}{\item \hskip-1.0cm}
\begin{list}{\ }{\setlength\leftmargin{1.0cm}}

\bibu \textsc{Antoniak, C.E.} (1974) Mixtures of Dirichlet processes with applications to Bayesian nonparametric problems. {\it Ann. Statist.}, 2, 1152--1174.

\bibu \textsc{Boender, C.G.E., Rinnooy Kan, A.H.G.} (1987) A multinomial Bayesian approach to the estimation of population and vocabulary size. {\it Biometrika}, 74, 849-856.

\bibu \textsc {Cerquetti, A.} (2011) Conditional $\alpha$-diversity for exchangeable Gibbs partitions driven by the stable subordinator. {\it Proceeding of the 7th Conference on Statistical Computation and Complex Systems, Padova, Italy, 2011}

\bibu \textsc {Chao, A., Bunge, J.} (2002) Estimating the number of species in a stochastic abundance model. {\it Biometrics}, 58, 531--539.

\bibu \textsc {Chao, A., Lee, S.M.} (1992) Estimating the number of classes via sample coverage. {\it J. Am. Statist. Assoc.} 87, 210-217. 

\bibu \textsc {Chao, A. Shen, T.} (2003) Nonparametric estimation of Shannon's index of diversity when there are unseen species in sample. {\it Envirom. Ecolog. Statistics}, 10, 429--443.

\bibu \textsc{Ewens, W. J.} (1972)  The sampling theory of selectively neutral alleles. {\it Theor. Pop. Biol.}, 3, 87-112.

\bibu \textsc{Favaro, S., Lijoi, A., Mena, R.H. and Pr\"unster, I.} (2009) Bayesian non-parametric inference for species variety with a two-parameter Poisson-Dirichlet process prior. {\it JRSS-B}, 71, 993-1008.

\bibu \textsc{Favaro, S., Lijoi, A. and Pr\"unster, I.} (2011) Asymptotics for a Bayesian nonparametric estimator of species variety. {\it Bernoulli} (to appear).

\bibu \textsc {Favaro, S., Lijoi, A. and Pr\"unster, I.} (2012) On the stick breaking representation of normalized inverse Gaussian priors. {\it Biometrika} (to appear)

\bibu \textsc{Ferguson, T.S.} (1973) A Bayesian analysis of some nonparametric problems. {\it Ann. Statist.}, 1, 209--230.

\bibu \textsc{Fisher, R.A., Corbet, A.S. and Williams, C. B.} (1943) The relation between the number of species and the number of individuals in a random sample of an animal population. {\it J. Animal Ecol.} 12, 42--58.

\bibu \textsc{Gill, C.A. and Joanes, D. N.} (1979) Bayesian estimation of Shannon's index of diversity. Biometrika, 66, 1, 81-85.


\bibu \textsc{Gnedin, A. and Pitman, J. } (2006) {Exchangeable Gibbs partitions  and Stirling triangles.} {\it Journal of Mathematical Sciences}, 138, 3, 5674--5685.

\bibu \textsc{Hill, B.M.} (1979) Posterior moments of the number of species in a finite population and the posterior probability of finding a new species. {\it J. Am. Statist. Assoc.}, 74, 668--673.

\bibu \textsc{James, L.F, Lijoi, A. and Pr\"unster, I.} (2009) Posterior analysis for normalized random measures with independent increments. {\it Scand. J. Statist.}, 36, 76--97.

\bibu \textsc{Lijoi, A., Mena, R.H. and Pr\"unster, I.} (2007) Bayesian nonparametric estimation of the probability of discovering new species.  {\it Biometrika}, 94, 769--786.

\bibu \textsc{Lijoi, A., Pr\"unster, I. and Walker, S.G.} (2008) Bayesian nonparametric estimator derived from conditional Gibbs structures. {\it Annals of Applied Probability}, 18, 1519--1547.

\bibu \textsc {Lloyd, M and Ghelardi, R. J.} (1964) A table for calculating the equitability component of species diversity. {\it J. Animal Ecol.} 33, 217--225.

\bibu \textsc{Pielou, E.C.} (1975) {\it Ecological Diversity} New York: Wiley.

\bibu \textsc{Pitman, J.} (1995) Exchangeable and partially exchangeable random partitions. {\it Probab. Th. Rel. Fields}, 102: 145-158.

\bibu \textsc{Pitman, J.} (1996) Some developments of the Blackwell-MacQueen urn scheme. In T.S. Ferguson, Shapley L.S., and MacQueen J.B., editors, {\it Statistics, Probability and Game Theory}, volume 30 of {\it IMS Lecture Notes-Monograph Series}, pages 245--267. Institute of Mathematical Statistics, Hayward, CA.

\bibu \textsc{Pitman, J.} (2003) {Poisson-Kingman partitions}. In D.R. Goldstein, editor, {\it Science and Statistics: A Festschrift for Terry Speed}, volume 40 of Lecture Notes-Monograph Series, pages 1--34. Institute of Mathematical Statistics, Hayward, California.

\bibu \textsc{Pitman, J.} (2006) {\it Combinatorial Stochastic Processes}. Ecole d'Et\'e de Probabilit\'e de Saint-Flour XXXII - 2002. Lecture Notes in Mathematics N. 1875, Springer.

\bibu \textsc{Pitman, J. and Yor, M.} (1997) The two-parameter Poisson-Dirichlet distribution derived from a stable subordinator. {\it Ann. Probab.}, 25, 855--900.

\bibu \textsc {Shen, T.J., Chao, A. and Lin, C.F.} (2003) Predicting the number of new species in further taxonomic sampling. {\it Ecology}, 84, 798-804.

\bibu \textsc {Simpson, E.H.} (1949) Measurement of diversity. {\it Nature} 163, 688

\end{list}

\end{document}